# The Erdös-Straus Conjecture and Pythagorean Primes

Bernd R. Schuh

Dr. Bernd Schuh, D-50968 Cologne, Germany, bernd.schuh@netcologne.de

**Abstract**

The Diophantine equation *4/n = 1/x + 1/y + 1/z* for a Pythagorean prime n is split into two independent Diophantine equations, which correspond to two different types of solution. The solvability of these equations forces certain restrictions on allowed Pythagorean primes. Empirical evidence suggests that these restrictions hold for all Pythagorean primes, which I state as two independent conjectures. One can be formulated as follows: every Pythagorean prime can be written as *p =(4ab-1)(4c-1) – 4ab²/ d*, where *a, b, c* are natural numbers and *d* is a divisor of *ab*. The second conjecture reads: every Pythagorean prime can be written as *p =(4ab-1)(4c-1) – 4ac*, where *a, b, c* are natural numbers. I give a new straightforward plausibility for the latter conjecture (which has been formulated independently by other authors) and I outline a practicle and effective algorithm to determine *a,b,c* for a given *p*.

**I Introduction**

In the following I will call an equation $\Phi(n,a,b,...,x,y,z)=0$ solvable in $n$ iff for a given $n \in \mathbb{N}$ positive integers $a,b,...,z$ exist which fulfill this equation. These integers are called the solution. A well known conjecture by Erdös-Straus (ESC in the following) states that the Diophantine equation

$$4/n - (1/x + 1/y + 1/z) = 0 \qquad (1)$$

is solvable in every integer $n \geq 2$. There is an impressive body of evidence for the validity of the conjecture, see e.g. [1,2,3,4], but no valid proof. Most approaches sought parametric solutions and covering of integers by arithmetic progression. Probably, different methods are needed. A key element appear to be Pythagorean primes (Pp in the following), since solving (1) is a nontrivial task only if $n$ is a Pythagorean prime. I list several known facts to support this assertion.

(i) If equ. (1) is solvable in $n$ it is also solvable in any multiple $n' = kn$ of $n$. Proof: simply set $x' = kx$ etc..

(ii) Equ. (1) is solvable in all even $n = 2k$. Proof: simply put $x = k, y = 2k = z$.

(iii) Equ. (1) is solvable in all $n$ of the form $n = 4k - 1$. Proof: set $x = kn, y = k(k+1), z = k+1$. (See [5] for similar results.)

These three facts constitute the claim that to prove the ESC it suffices to show that a solution exists for all primes of the form $n = 1 + 4K$, i.e. all Pp. We neglect further progress on restrictions on $K$, like the well known ones by Mordell [6] and others [7], since they are of no importance for the following arguments. A fourth fact, however, is notable and important in the following:

(iv) There are exactly two types of solutions, say (A) and (B). In type (A) only one of the variables $x, y, z$ is proportional to $n$. In type (B) two variables are multiples of $n$. [2]

**II Solutions**

In this section I will deduce conditions that guaranty solvability of equation (1) in Pythagorean primes. For the rest of the paper I use the abbreviation $\beta_w := 4w - 1$ for elements of the residue class $3 \bmod 4$. The restriction to Pythagorean primes leads to more transparent results than earlier approaches using similar reformulations of equation (1), see e.g. [9,10]. Central to our arguments is a decomposition of (1) into two equations which correspond to one of the solution types (A) or (B) mentioned in fact (iv) and introduced by Bernstein [2].

*Theorem 1*

*The diophantine equation (1) is solvable in $n$ being a prime if and only if either*

$$\Phi_A(n,a,d,z) := \beta_a z - n(d+a) = 0 \tag{2a}$$

*or*

$$\Phi_B(n,a,d,z) := \beta_a z - d - na = 0 \tag{2b}$$

*is solvable, where $d$ is a divisor of $a^2$. The solution of (1) then reads*

$$z, \; y = az/d, \; x = an \tag{3a}$$

*for (2a), and*

$$z, \; y = naz/d, \; x = an \tag{3b}$$

*for (2b).*

*Proof*: In a first step we show that (1) is equivalent to equations (2a) and (2b) for primes. According to fact (iv) at least one of the variables $x, y, z$ is proportional to $n$. Thus we write $x = an$, then

$$4/n = 1/x + 1/y + 1/z \Leftrightarrow$$
$$\beta_a yz = na(y+z) \Leftrightarrow$$
$$\beta_a^2 yz - \beta_a na(y+z) + n^2 a^2 = n^2 a^2 \Leftrightarrow$$
$$(\beta_a z - na)(\beta_a y - na) = n^2 a^2$$

Since $n$ is a prime either $n^2$ divides one of the factors on the l.h.s. of the last equation or $n$ divides both. In the later case we have $\beta_a z - na = nd$ and $\beta_a y - na = nd_c$ and thus (2a) and (3a), in the former case $\beta_a z - na = d$ and $\beta_a y - na = n^2 d_c$ and thus (2b), (3b). In both cases $d$ denotes a divisor of $a^2$ and $d_c$ its cofactor: $a^2 = dd_c$.

It remains to be shown that $az/d$ in the formula for $y$ is a positive integer. Since any divisor of $a^2$ may be written as a product of divisors of $a$, $d = a_1 a_2$, one can write $\beta_a z/a_1 = n(a_2 + a_{1c})$ instead of (2a) and $\beta_a z/a_1 = a_2 + na_{1c}$ instead of (2b), where the subscript c denotes cofactors, as before. Since $(\beta_a, a_1) = 1$, one has $z/a_1 \in \mathbb{N}$ and thus $az/d = az/a_1 a_2 = a_{2c}(z/a_1) \in \mathbb{N}$.

In effect theorem 1 states that finding positive integers $x, y, z$ which fulfill (1) for primes $n$ is equivalent to finding positive integers $z, a, d$ which fulfill one of the equ. (2a) or (2b). The solutions of the two approaches are then related by equs. (3).

A consequence of equations (3a), (3b) is the

*Corollary:*

If $n$ is a Pythagorean prime, i.e. $n = 1 + 4K$, every solution of (3a) is a solution of type (A) in Bernsteins notation [2] (see fact (iv) in the introduction).

Every solution of (3b) is a solution of type (B).

proof: the second statement is obvious from (3b). To prove the first part of the corollary we observe that according to (2a) the prime $n$ divides either $\beta_a$ or $z$. The latter possibility can be excluded, because in that case, according to (3a), all three variables would be proportional to $n$, which contradicts equ. (1). Thus $n$ divides $\beta_a$. Since $n \equiv 1 \bmod 4$ the proportionality factor must be congruent $-1 \bmod 4$. So a parameter $b$ exists with

$$\beta_a = n\beta_b \Leftrightarrow$$
$$a = bn - K \Leftrightarrow \qquad (4)$$
$$a = b + \beta_b K$$

Therefore $(a, n) = 1$ and $y$ in equ. (3a) cannot have a divisor $n$.

The problem of solving equ. (1) for a given prime $n$ thus amounts to finding parameters $a$ and a divisor $d$ of $a^2$ which allow for positive integer solutions $z$ either of equ. (2a) or (2b). This reformulation of the ESC problem leads to interesting parametrizations of allowed primes. We state this fact as two separate theorems, depending on whether (2a) or (2b) holds:

*Theorem 2A*

*Equation (2a) is solvable in a prime $n = 1 + 4K$ if $K$ can be written*

$$K = \mu(\kappa\beta_b - 1) - \kappa b \qquad (5)$$

*with parameters $b, \mu, \kappa \in \mathbb{N}$.*

*Proof of theorem 2A*: We assume that $b, \mu, \kappa$ and $K$ are given as in (5). We <u>define</u> $a$ via (4) and get the factorization: $a = (\mu\beta_b - b)(\kappa\beta_b - 1)$. Defining furthermore $d := (\mu\beta_b - b)$, which obviously is a divisor of $a^2$, one has $a + d = (\mu\beta_b - b)\kappa\beta_b$. Since (2a) yields

$$\beta_b z = a + d, \qquad (6)$$

setting $z := \kappa(\mu\beta_b - b)$ completes the proof.

Unfortunately, condition (5) is only sufficient to solve (2a) but not necessary. To see this we assume (2a) to hold with $z, a, d$ given. We have already proven that a positive integer $b$ exists, connected to $a$ via (4), such that (2a) is equivalent to (6). Next we observe

*Lemma 1*: $a$ must be composite.

*Proof*: If it were a prime then $d \in \{1, a, a^2\}$. For $d = 1$ or $d = a = b + \beta_b K$ equation (6) yields either $b + 1 = \nu\beta_b$ or $2b = \nu\beta_b$ for some $\nu \in \mathbb{N}$. Both cannot be fulfilled with $b \in \mathbb{N}$. For $d = a^2$ we have from (4) $a^2 \equiv b^2 \mod \beta_b$ and one gets from (6) $\nu\beta_b - b = b^2$ with some $\nu \in \mathbb{N}$. But $4\nu b - \nu - b$ never yields a perfect square (see e.g. [8,11]). Thus $a$ must be composite.

Since from (4) $a \equiv b \mod \beta_b$ and $(b, \beta_b) = 1$ one can make a general Ansatz

$$a = (\mu\beta_b - b_1)(\kappa\beta_b - b_2) \qquad (7)$$

with integer parameters $\mu, \kappa \in \mathbb{N}$ and $b_1, b_2 \in \mathbb{Z}$. The latter must fulfill $b_1 b_2 = b + \lambda\beta_b$ according to (4). Solving (4) for $K$ with $a$ replaced by (7) one gets

$$K = \mu\kappa\beta_b - b_1\kappa - b_2\mu + \lambda \qquad (8)$$

as a necessary condition for the solvability of (2a). Of course, a six parameter-family of solutions is not an achievement. But the foregoing steps help to understand the idea behind (5), which one gets from (8) by the choice $b_1 = b, b_2 = 1, \lambda = 0$.

The situation for equ. (2b) is different. Here we are able to identify a four-parameter-family of solutions which is both sufficient and necessary for the solvability of (2b).

*Theorem 2B*

*Equation (2b) is solvable in a prime $n = 1 + 4K$ if and only if $K$ can be written as*

$$K = \mu \beta_a - a - d \qquad (9)$$

*with parameters $a, \mu \in \mathbb{N}$ and a divisor $d$ of $a^2$.*

*Proof of theorem 2B*: consider equ. (2b) and assume that a positive integer solution exists for a given $n = 1 + 4K$. Inserting this in $an + d$ one gets from (2b) $\beta_a(z - K) = d + a + K$. With a new parameter $\mu := z - K$ one can write $K = \mu \beta_a - a - d$, as claimed. If vice versa $n = 1 + 4K$ is given by equ. (9) in positive integers $a, \mu, d$ one inserts this expression on the r.h.s. of (2b) and defines $z$ by $z := K + \mu$, which is a positive integer and fulfills (2b). Done.

Note that the proof of theorem 2B does not necessitate $n$ to be a prime. It is, however, necessary for (2b) being equivalent to (1) (together with (2a)).

**III Discussion**

Up to now I have identified a three-parameter family of solutions for equation (2a) and a four-parameter family of solutions for (2b), the latter being not only sufficient but also necessary. The first always leads to solutions of type (A), i.e. with exactly one of the variables in (1) being proportional to $n$, whereas the second produces solutions of type (B) with two variables being multiples of $n$.

What is the further relevance of these parametrizations? Let us consider theorem 2A first. We observe that equation (5) is not trivial in the sense that it comprizes *all* odd numbers of the form $1 + 4K$. There are exceptions, e.g. all perfect squares and some numbers with $K \equiv 2 \mod 3$ like $K = 26$, and $K = 200$.

Secondly, we observe that it is usable in that it gives simple deductions for well known facts. E.g. it shows that equation (1) is solvable for all odd $K$. To see this, set $\mu = b = 1$ in (5) to get $K = 2\kappa - 1$. Thus theorem 2A ensures that whenever $K$ is an odd integer the Diophantine equation is solvable. It is also easy to see that $K \equiv 1 \mod 3$ is easily represented by (5). Set $\mu = 1 = \kappa$ to get $K = 3b - 2$.

Furthermore, $K \equiv 2 \bmod 3$ leads to $n \equiv 0 \bmod 3$ which is not a prime. Thus only $K \equiv 0 \bmod 3$ and $K \equiv 0 \bmod 2$ remain as interesting possibilities, a long known restriction, which states that exceptions (if any) to the ESC must be sought among the primes $n \equiv 1 \bmod 24$.

Surprisingly, the parametrization (5) worked for all Pythagorean primes we tried, extremely well. We were able to calcultate $b$, $\mu$ and $K$ even for many-digit primes with a pocket calculator (and open source prime factorization programs). The reason is that the derivation of (5) offers a convenient algorithm to determine the parameters in this equation for a given $K$. It works as follows:

*Algorithm to determine $b$, $\mu$ and $K$*

For given $n$ or $K$ respectively calculate $a$ from equ. (4), starting conveniently with $b=1$.
Then factorize $a$ and check whether one of the factors is congruent $-1 \bmod \beta_b$. If so, the cofactor will be congruent $-b \bmod \beta_b$ and from (8) with the special choice $b_1 = b, b_2 = 1, \lambda = 0$ one reads off $\mu$ and $K$. Store the solution $b, \mu, K$.
Repeat the procedure with $b \rightarrow b+1$.
The process ends, when $b+1$ exceeds a limit, e.g. $\lceil (2+K)/3 \rceil$. Then all solutions are in store. If the store is empty there is no solution and $K$ cannot be represented by (5).

The limit can be derived from a rough estimation: First rewrite (5) as
$$n = \kappa \beta_b \beta_\mu - \beta_\mu - \kappa \qquad (10)$$
Then $n > \kappa(\beta_\mu(\beta_b - 1) - 1)$. And $1 + n/\kappa > 3(\beta_b - 1)$ yields the cited limit.

Take $n = 560281$ as a nontrivial example illustrating the algorithm. With $b = 1$ one gets $3K + 1 = 420211 = 11 \times 38201$ and therefore $(1+11)/3 = 4 = \mu$ and $(1+38201)/3 = 12734 = \kappa$. $b = 2$ yields a second solution: $b = 2, \mu = 10, \kappa = 2060$. To find *all* solutions one would have to go to $b = 46690$.

Since (5) apparently represents the Pythagorean primes so extremely well, it is tempting to put forward the following conjecture:

*Conjecture A*

All primes of the form $1+4K$ belong to the set $S_A$ defined by
$$S_A = \{abc - b - c : a, b \equiv 3 \bmod 4; c \in \mathbb{N}\}$$
i.e. the parametrization given in theorem 2A.

Since any member of $S_A$ gives rise to a solution of equation (1) as stated in theorem 2A it is clear that the ESC is valid if the conjecture is.

Up to now, conjecture A is unproven. After all, one can show

*Lemma 2*: $S_A$ contains infinitely many primes.

*Proof*: We use a theorem by Iwaniec [12]. Consider (10) as a quadratic polynomial in two variables $x = \kappa$ and $y = \beta_\mu$. Then obviously $n(x,y) = \beta_b xy - x - y$ is irreducible and $\partial n / \partial x = \beta_b y - 1$ and $\partial n / \partial y = \beta_b x - 1$ are linearly independent. Thus the assumptions of theorem 1 in [12] are fulfilled. According to conclusion (i) of that theorem the number of primes up to a given $N$ being represented by $n(x,y)$ grows with $N$ faster than $N / \log N$.

After completion of this work I noticed that the parametrization (5) was put forward in [10] already. The authors of [10] made extensive calculations and verified conjecture A for Pythagorean primes up to $n < 10^{14}$. They have also proven that the set $S_A$ does not contain perfect squares, and they list a finite number of non-primes which are not members of $S_A$.

I would like to point out, however, that my "derivation" of the conjecture is more straightforward and leads to a usable algorithm to identify given Pythagorean primes as members of the three-parameter family of solutions.

Ad theorem 2B.

To illustrate the result of theorem 2B take e.g. $a = 1 = d$ in (9) and let $\mu$ run from 1 to 5. This generates the Pythagorean primes 5, 17, 29, 41, 53. The primes 13 and 37 are produced by the choice $d = a$, setting $\mu = 1$ and letting $a$ run.

Depending on the value of $a$ there can be a huge amount of choices for the parameter $d$. Writing quite generally $d = a_1 a_2$, we have $a = \lambda a_1 = \nu a_2$ since $a_{1,2}$ are divisors of $a$ and condition (9) becomes

$$n = \beta_{\lambda a_1} \beta_\mu - 4\lambda a_1^2 / \nu \qquad (11)$$

While $\lambda$, $a_1$ and $\mu$ are independent positive integers, $\nu$ is restricted to divisors of $\lambda a_1$. With this restriction in mind, the existence of parameters $\lambda, a_1, \mu, \nu$ relating to $n$ via (11) is a necessary <u>and</u> sufficient condition for the solvability of equation (2b).

Other than theorem 2A, theorem 2B and (11) is of restricted usefulness to search for solutions of equ. (1) with $n$ given, since it offers no obvious algorithm for limiting the variety of choices. Nonetheless, empirically the parametrization (11) works extremely well. All Pythagorean primes we checked had a representation in terms of equ. (11). Thus we propose a second conjecture based on equ. (2b) and theorem 2B:

*Conjecture B*

All primes of the form $1+4K$ belong to the set $S_B$ defined by

$$S_B = \left\{ \beta_{\lambda\alpha}\beta_{\lambda\mu} - \frac{4\lambda\alpha^2}{\nu} : \lambda, \alpha, \mu, \nu, \frac{\lambda\alpha}{\nu} \in \mathbb{N} \right\}$$

i.e. the parametrization given in equ. (11).

Also conjecture B is not trivial. The set $S_B$ does not include <u>all</u> odd numbers congruent 1mod4, e.g. it does not include perfect squares, too.

**Conclusion**.

In a first step we have replaced the Diophantine equation (1) by two simpler equations (2a) and (2b) corresponding to the two types of solutions known as type (A) and type (B) [2]. In type (A) one of the variables is proportional to $n$, in type (B) two of them are. Each equation leads to a different parametric representation of Pythagorean primes, $p \in S_A$ and $p \in S_B$, see equations (10) and (11), respectively. Both are sufficient to guaranty the solvability of the corresponding equation. But only the type (B) conjecture is also necessary. There is rather strong evidence that each of these representations comprizes <u>all</u> Pythagorean primes. We thus propose two corresponding conjectures A and B. Should both conjectures turn out to be valid the Erdös-Straus conjecture were proved. Should the ESC be proven then at least conjecture B is bound to hold.